# LARGE DEVIATIONS FOR RENORMALIZED SELF-INTERSECTION LOCAL TIMES OF STABLE PROCESSES

By Richard Bass[1], Xia Chen[2] and Jay Rosen[3]

*University of Connecticut, University of Tennessee and College of Staten Island, CUNY*

We study large deviations for the renormalized self-intersection local time of $d$-dimensional stable processes of index $\beta \in (2d/3, d]$. We find a difference between the upper and lower tail. In addition, we find that the behavior of the lower tail depends critically on whether $\beta < d$ or $\beta = d$.

**1. Introduction.** Let $X_t$ be a nondegenerate $d$-dimensional stable process of index $\beta$. We assume that $X_t$ is symmetric, that is, $X_t \stackrel{d}{=} -X_t$, but we do not assume it is spherically symmetric. Thus,

$$E(e^{i\lambda \cdot X_t}) = e^{-t\psi(\lambda)}, \tag{1.1}$$

where $\psi(\lambda) \ge 0$ is continuous, positively homogeneous of degree $\beta$, that is, $\psi(r\lambda) = r^\beta \psi(\lambda)$ for each $r \ge 0$, $\psi(-\lambda) = \psi(\lambda)$ and for some $0 < c < C < \infty$,

$$c|\lambda|^\beta \le \psi(\lambda) \le C|\lambda|^\beta. \tag{1.2}$$

In studying the self intersections of $\{X_t; t \ge 0\}$, one is naturally led to try to give meaning to the formal expression

$$\int_0^t \int_0^s \delta_0(X_s - X_r)\, dr\, ds, \tag{1.3}$$

where $\delta_0(x)$ is the Dirac delta "function." Let $\{f_\varepsilon(x); \varepsilon > 0\}$ be an approximate identity and set

$$\int_0^t \int_0^s f_\varepsilon(X_s - X_r)\, dr\, ds. \tag{1.4}$$

---

Received October 2003; revised May 2004.
[1]Supported in part by NSF Grant DMS-02-44737.
[2]Supported in part by NSF Grant DMS-01-02238.
[3]Supported in part by grants from the NSF and from PSC-CUNY.

When $\beta > d$, so that necessarily $d = 1$ and $\{X_t; t \geq 0\}$ has local times $\{L_t^x; (x,t) \in R^1 \times R^1_+\}$, (1.4) converges as $\varepsilon \to 0$ to $\frac{1}{2}\int (L_t^x)^2\,dx$. Large deviations for this object have been studied in [7].

In this paper we assume that $\beta \leq d$. In this case (1.4) blows up as $\varepsilon \to 0$. We consider instead

$$\gamma_{t,\varepsilon} = \int_0^t \int_0^s f_\varepsilon(X_s - X_r)\,dr\,ds - E\left\{\int_0^t \int_0^s f_\varepsilon(X_s - X_r)\,dr\,ds\right\} \tag{1.5}$$

and let

$$\gamma_t = \lim_{\varepsilon \to 0} \gamma_{t,\varepsilon} \tag{1.6}$$

whenever the limit exists. It is known that this happens if (and only if) $\beta > 2d/3$, and then $\gamma_t$ is continuous in $t$ almost surely [22, 23, 26]. In this case we refer to $\gamma_t$ as the renormalized self-intersection local time for the process $X_t$. Renormalized self-intersection local time, originally studied by Varadhan [28] for its role in quantum field theory, turns out to be the right tool for the solution of certain "classical" problems such as the asymptotic expansion of the area of the Wiener and stable sausages in the plane and fluctuations of the range of stable random walks. See [14, 15, 18, 25]. In [27] we show that $\gamma_t$ can be characterized as the continuous process of zero quadratic variation in the decomposition of a natural Dirichlet process. For further work on renormalized self-intersection local times, see [3, 10, 16, 21, 26].

The goal of this paper is to study the large deviations of $\gamma_t$, generalizing the recent work for planar Brownian motion of the first two authors [2].

THEOREM 1. *Let $X_t$ be a symmetric stable process of order $2d/3 < \beta \leq d$ in $R^d$. Then, for some $0 < a_\psi < \infty$ and any $h > 0$,*

$$\lim_{t\to\infty} \frac{1}{t} \log P(\gamma_t \geq ht^2) = -h^{\beta/d} a_\psi. \tag{1.7}$$

The constant $a_\psi$ is described in Section 4 and is related to the best possible constant in a Gagliardo–Nirenberg type inequality.

$\gamma_t$ is not symmetric. In fact, the lower tail has very different behavior.

THEOREM 2. *Let $X_t$ be a symmetric stable process of order $\beta > 2d/3$ in $R^d$. Then we can find some $0 < b_\psi < \infty$ such that if $\beta < d$,*

$$\lim_{t\to\infty} \frac{1}{t} \log P(-\gamma_t \geq t) = -b_\psi, \tag{1.8}$$

*while if $\beta = d$,*

$$\lim_{t\to\infty} \frac{1}{t} \log P(-\gamma_1 \geq p_1(0) \log t) = -b_\psi, \tag{1.9}$$

*where $p_t(x)$ is the continuous density function for $X_t$.*

We are unable to identify the constant $0 < b_\psi < \infty$.

Using the scaling property $\{X(ts); s \geq 0\} \stackrel{d}{=} t^{1/\beta}\{X(s); s \geq 0\}$ of the stable process, it is easy to check that

$$\gamma_t \stackrel{d}{=} t^{2-d/\beta}\gamma_1. \tag{1.10}$$

Thus, (1.7)–(1.9) are equivalent to

$$\lim_{h\to\infty} \frac{1}{h^{\beta/d}} \log P(\gamma_1 \geq h) = -a_\psi, \tag{1.11}$$

$$\lim_{h\to\infty} \frac{1}{h^{\beta/(d-\beta)}} \log P(-\gamma_1 \geq h) = -b_\psi, \qquad \beta \in (2d/3, d), \tag{1.12}$$

$$\lim_{h\to\infty} \frac{1}{e^{p_1(0)h}} \log P(-\gamma_1 \geq h) = -b_\psi, \qquad \beta = d. \tag{1.13}$$

Equations (1.11) and (1.12) show that

$$\lim_{h\to\infty} \frac{1}{h} \log P(|\gamma_1|^{\beta/d} \geq h) = -a_\psi, \tag{1.14}$$

which implies that

$$E(e^{\lambda|\gamma_1|^{\beta/d}}) \begin{cases} < \infty, & \text{if } \lambda < a_\psi^{-1}, \\ = \infty, & \text{if } \lambda > a_\psi^{-1}. \end{cases} \tag{1.15}$$

Our large deviation results lead to the following law of the iterated logarithm (LIL) type results.

THEOREM 3. *Let $X_t$ be a symmetric stable process of order $2d/3 < \beta \leq d$ in $R^d$. Then*

$$\limsup_{t\to\infty} \frac{\gamma_t}{t^{(2-d/\beta)}(\log\log t)^{d/\beta}} = a_\psi^{-d/\beta} \qquad a.s. \tag{1.16}$$

THEOREM 4. *Let $X_t$ be a symmetric stable process of order $\beta > 2d/3$ in $R^d$. If $\beta < d$, then*

$$\liminf_{t\to\infty} \frac{\gamma_t}{t^{(2-d/\beta)}(\log\log t)^{d/\beta-1}} = -b_\psi^{-(d/\beta-1)} \qquad a.s., \tag{1.17}$$

*while if $\beta = d$, then*

$$\liminf_{t\to\infty} \frac{1}{t\log\log\log t}\gamma_t = -p_1(0) \qquad a.s. \tag{1.18}$$

The methods needed for this paper are very different from those used in [2] for planar Brownian motion. In that case, and more generally when $\beta = d$, the upper bound for large deviations for $\gamma_t$ comes from a soft argument involving scaling. This argument breaks down when $\beta < d$. Instead, we obtain the upper bound using careful moment arguments developed in Sections 2 and 3.

Another major difference between this paper and [2] is in the proof of the lower bound for large deviations for $-\gamma_t$ when $\beta < d$. Suppose we divide the time interval $[0, n]$ into subintervals $I_k = [k, k+1]$, $k = 0, \dots, n-1$, let $\Gamma(I_k)$ denote renormalized self-intersection local time for the piece of the path generated by times in $I_k$, and let $A(I_j; I_k)$ denote the intersection local time for the two pieces generated by times in $I_j$ and $I_k$ when $j \neq k$. Then the major contribution to the renormalized self-intersection local time for planar Brownian motion on the interval $[0, n]$ comes from $\sum_{j<k}[A(I_j; I_k) - EA(I_j; I_k)]$; the contribution from $\sum_k \Gamma(I_k)$ is smaller. In contrast, when $\beta < d$, both contributions are of the same order of magnitude. As a result, the lower bound for $-\gamma_t$ when $\beta < d$ requires a much more delicate argument.

Our paper is organized as follows. In Section 2 we obtain bounds on exponential moments of the intersection local time for two independent processes, which is then used in Section 3, following an approach due to Le Gall, to obtain bounds on exponential moments of the renormalized self-intersection local time $\gamma_t$, and, in particular, to obtain an exponential approximation of $\gamma_t$ by its regularization $\gamma_{t,\varepsilon}$. Together with some results from [8], this allows us to prove Theorem 1 in Section 4. In Sections 5 and 6 we prove Theorem 2 on the lower tail of $\gamma_t$. Finally, these results are used in Sections 7 and 8 to prove the LILs of Theorems 3 and 4, respectively.

## 2. Intersection local times.

Let $X_t, X'_t$ be two independent copies of the symmetric stable process of order $\beta$ in $R^d$ with characteristic exponent $\psi$ and set

$$\alpha_{t,\varepsilon} \stackrel{\text{def}}{=} \int_0^t \int_0^t \int_{R^d} f_\varepsilon(X_s - X'_r)\, dr\, ds, \tag{2.1}$$

where $f_\varepsilon$ is an approximate $\delta$—function at zero, that is, $f_\varepsilon(x) = f(x/\varepsilon)/\varepsilon^d$ with $f \in \mathcal{S}(R^d)$ a positive, symmetric function with $\int f\, dx = 1$. If $\widehat{f}(p)$ denotes the Fourier transform of $f$, then $\widehat{f}(\varepsilon p)$ is the Fourier transform of $f_\varepsilon$ and we have, from (2.1),

$$\alpha_{t,\varepsilon} = (2\pi)^{-d} \int_0^t \int_0^t \int_{R^d} e^{ip\cdot(X_s - X'_r)} \widehat{f}(\varepsilon p)\, dp\, dr\, ds. \tag{2.2}$$

THEOREM 5. *Let $X_t, X'_t$ be independent copies of a symmetric stable process of order $d/2 < \beta \le d$ in $R^d$. Then for all $\rho > 0$ sufficiently small, we can find some $\theta > 0$ such that*

$$\sup_{\varepsilon, \varepsilon', t > 0} E\left(\exp\left\{\theta \left|\frac{\alpha_{t,\varepsilon} - \alpha_{t,\varepsilon'}}{|\varepsilon - \varepsilon'|^\rho t^{2-(d+\rho)/\beta}}\right|^{\beta/(d+\rho)}\right\}\right) < \infty. \tag{2.3}$$

*Furthermore,*

$$\lim_{\theta\to 0}\sup_{\varepsilon,\varepsilon',t>0} E\biggl(\exp\biggl\{\theta\biggl|\frac{\alpha_{t,\varepsilon}-\alpha_{t,\varepsilon'}}{|\varepsilon-\varepsilon'|^{\rho}t^{2-(d+\rho)/\beta}}\biggr|^{\beta/(d+\rho)}\biggr\}\biggr)=1. \tag{2.4}$$

PROOF. From (2.2), we have that

$$\alpha_{t,\varepsilon}-\alpha_{t,\varepsilon'}=(2\pi)^{-d}\int_0^t\int_0^t\int_{R^d} e^{ip\cdot(X_s-X'_r)}(\widehat{f}(\varepsilon p)-\widehat{f}(\varepsilon' p))\,dp\,dr\,ds. \tag{2.5}$$

Hence,

$$\begin{aligned} &E(\{\alpha_{t,\varepsilon}-\alpha_{t,\varepsilon'}\}^n)\\ &\quad=(2\pi)^{-nd}\int_{[0,t]^n}\int_{[0,t]^n}\int_{R^{dn}} E(e^{i\sum_{k=1}^n p_k(X_{s_k}-X'_{r_k})})\\ &\qquad\qquad\times\prod_{j=1}^n\{\widehat{f}(\varepsilon p_j)-\widehat{f}(\varepsilon' p_j)\}\,dp_j\,dr_j\,ds_j. \end{aligned} \tag{2.6}$$

We then use the decomposition

$$[0,t]^n\times[0,t]^n=\bigcup_{\pi,\pi'}D_n(\pi,\pi'),$$

where the union runs over all pairs of permutations $\pi,\pi'$ of $\{1,\dots,n\}$ and $D_n(\pi,\pi')=\{(r_1,\dots,r_n,s_1,\dots,s_n)|r_{\pi_1}<\cdots<r_{\pi_n}\le t,\, s_{\pi'_1}<\cdots<s_{\pi'_n}\le t\}$. Using this, we then obtain

$$\begin{aligned} &E(\{\alpha_{t,\varepsilon}-\alpha_{t,\varepsilon'}\}^n)\\ &\quad=(2\pi)^{-nd}\sum_{\pi,\pi'}\int_{D_n(\pi,\pi')}\int_{R^{dn}} E(e^{i\sum_{k=1}^n p_k(X_{s_k}-X'_{r_k})})\\ &\qquad\qquad\times\prod_{j=1}^n\{\widehat{f}(\varepsilon p_j)-\widehat{f}(\varepsilon' p_j)\}\,dp_j\,dr_j\,ds_j. \end{aligned} \tag{2.7}$$

On $D_n(\pi,\pi')$, we can write

$$\begin{aligned} &\sum_{k=1}^n p_k(X_{s_k}-X'_{r_k})\\ &\quad=\sum_{k=1}^n u_{\pi,k}(X_{r_{\pi_k}}-X_{r_{\pi_{k-1}}})-\sum_{k=1}^n u_{\pi',k}(X'_{s_{\pi'_k}}-X'_{s_{\pi'_{k-1}}}), \end{aligned} \tag{2.8}$$

where, for any permutation $\pi$, we set $u_{\pi,k}=\sum_{j=k}^n p_{\pi_j}$. Hence, on $D_n(\pi,\pi')$,

$$\begin{aligned} &E(e^{i\sum_{k=1}^n p_k(X_{s_k}-X'_{r_k})})\\ &\quad=e^{-\sum_{k=1}^n\psi(u_{\pi,k})(r_{\pi_k}-r_{\pi_{k-1}})}e^{-\sum_{k=1}^n\psi(u_{\pi',k})(s_{\pi'_k}-s_{\pi'_{k-1}})}. \end{aligned} \tag{2.9}$$

We will use the bound $|\widehat{f}(\varepsilon p_j)-\widehat{f}(\varepsilon' p_j)|\le C|\varepsilon-\varepsilon'|^\rho|p_j|^\rho$ for any $\rho\le 1$. Using the Cauchy–Schwarz inequality, we have

$$\begin{aligned}
&\int_{R^{dn}} E(e^{i\sum_{k=1}^n p_k(X_{s_k}-X'_{r_k})})\prod_{j=1}^n |p_j|^\rho\,dp_j\\
&\qquad\le\left(\int_{R^{dn}} e^{-2\sum_{k=1}^n \psi(u_{\pi,k})(r_{\pi_k}-r_{\pi_{k-1}})}\prod_{j=1}^n|p_j|^\rho\,dp_j\right)^{1/2}\\
&\qquad\quad\times\left(\int_{R^{dn}} e^{-2\sum_{k=1}^n \psi(u_{\pi',k})(s_{\pi'_k}-s_{\pi'_{k-1}})}\prod_{j=1}^n|p_j|^\rho\,dp_j\right)^{1/2}.
\end{aligned}\tag{2.10}$$

Now $\prod_{j=1}^n|p_j|=\prod_{j=1}^n|p_{\pi_j}|=\prod_{j=1}^n|u_{\pi,j}-u_{\pi,j+1}|\le\prod_{j=1}^n|u_{\pi,j}|+|u_{\pi,j+1}|$ so that, using (1.2) for the second inequality,

$$\begin{aligned}
&\int_{R^{2n}} e^{-2\sum_{k=1}^n \psi(u_{\pi,k})(r_{\pi_k}-r_{\pi_{k-1}})}\prod_{j=1}^n|p_j|^\rho\,dp_j\\
&\qquad\le\sum_h\int_{R^n} e^{-2\sum_{k=1}^n \psi(u_{\pi,k})(r_{\pi_k}-r_{\pi_{k-1}})}\prod_{j=1}^n|u_{\pi,j}|^{h_j\rho}\,du_{\pi,j}\\
&\qquad\le\sum_h\int_{R^n} e^{-c\sum_{k=1}^n |u_{\pi,k}|^\beta(r_{\pi_k}-r_{\pi_{k-1}})}\prod_{j=1}^n|u_{\pi,j}|^{h_j\rho}\,du_{\pi,j}\\
&\qquad\le C^n\sum_h\prod_{j=1}^n(r_{\pi_k}-r_{\pi_{k-1}})^{-(d+h_j\rho)/\beta},
\end{aligned}\tag{2.11}$$

where the sum runs over all $h=(h_1,\dots,h_n)$ such that each $h_j=0,1$ or $2$ and $\sum_{j=1}^n h_j=n$.

Hence, taking $\rho>0$ sufficiently small that $(d+2\rho)/2\beta<1$, we have

$$\begin{aligned}
&E\left(\left|\frac{\alpha_{t,\varepsilon}-\alpha_{t,\varepsilon'}}{|\varepsilon-\varepsilon'|^\rho}\right|^n\right)\\
&\qquad\le C^n(n!)^2\left(\sum_h\int_{r_1<\cdots<r_n\le t}\prod_{j=1}^n(r_j-r_{j-1})^{-(d+h_j\rho)/2\beta}\,dr_j\right)^2\\
&\qquad\le C^n\left(t^{n(1-(d+\rho)/2\beta)}\frac{n!}{\Gamma(n(1-(d+\rho)/2\beta))}\right)^2\\
&\qquad\le C^n t^{2n(1-(d+\rho)/2\beta)}(n!)^{(d+\rho)/\beta}.
\end{aligned}\tag{2.12}$$

Hence, by Hölder's inequality,

$$\begin{aligned}
E\left(\left|\frac{\alpha_{t,\varepsilon}-\alpha_{t,\varepsilon'}}{|\varepsilon-\varepsilon'|^\rho t^{2-(d+\rho)/\beta}}\right|^{n\beta/(d+\rho)}\right)&\le E\left(\left|\frac{\alpha_{t,\varepsilon}-\alpha_{t,\varepsilon'}}{|\varepsilon-\varepsilon'|^\rho t^{2-(d+\rho)/\beta}}\right|^n\right)^{\beta/(d+\rho)}\\
&\le C^n n!.
\end{aligned}\tag{2.13}$$

Theorem 5 follows easily from this. □

If we set

$$\alpha_{s,t,\varepsilon} \stackrel{\text{def}}{=} \int_0^s \int_0^t f_\varepsilon(X_s - X'_r)\, dr\, ds, \tag{2.14}$$

then by the same method we can show that

$$\alpha_{s,t} = \lim_{\varepsilon \to 0} \alpha_{s,t,\varepsilon} \tag{2.15}$$

exists a.s. and in all $L^p$ spaces and for some $\theta > 0$,

$$\sup_{s,t>0} E\left(\exp\left\{\theta \left|\frac{\alpha_{s,t}}{(st)^{1-d/2\beta}}\right|^{\beta/d}\right\}\right) < \infty. \tag{2.16}$$

Let $p_t(x)$ denote the density function for $X_t$ started at the origin.

THEOREM 6. *Let $X_t, X'_t$ be independent copies of a symmetric stable process of order $d/2 < \beta < d$ in $R^d$. Let $P^{(x_0,y_0)}$ be the joint law of $(X_t, X'_t)$ when $X_t$ is started at $x_0$ and $X'_t$ is started at $y_0$. Then*

$$E^{(x_0,y_0)}(\alpha_{s,t}) \le c_\psi[s^{2-d/\beta} + t^{2-d/\beta} - (s+t)^{2-d/\beta}], \tag{2.17}$$

*where*

$$c_\psi = \frac{p_1(0)}{(d/\beta - 1)(2 - d/\beta)}. \tag{2.18}$$

*If $x_0 = y_0$, then we have equality in (2.17).*

*If $\beta = d$, then we obtain*

$$E^{(x_0,y_0)}(\alpha_{s,t}) \le p_1(0)[(s+t)\log(s+t) - t\log t - s\log s] \tag{2.19}$$

*with equality if $x_0 = y_0$.*

PROOF. We have

$$\begin{aligned} &E^{(x_0,y_0)}\left(\int_0^s \int_0^t f_\varepsilon(X_r - X'_u)\, dr\, du\right) \\ &\quad = \int_0^s \int_0^t \int f_\varepsilon(x-y) p_r(x-x_0) p_u(y-y_0)\, dx\, dy\, dr\, du \\ &\quad = \int_0^s \int_0^t \int f_\varepsilon(x) p_r(x+y-(x_0-y_0)) p_u(y)\, dx\, dy\, dr\, du \\ &\quad = \int_0^s \int_0^t \int f_\varepsilon(x) p_{r+u}(x-(x_0-y_0))\, dx\, dr\, du, \end{aligned} \tag{2.20}$$

where the last line follows from the semigroup property. Letting $\varepsilon \to 0$ and using the fact that (2.15) converges in $L^1$,

$$E^{(x_0,y_0)}(\alpha_{s,t}) = \int_0^s \int_0^t p_{r+u}(x_0 - y_0)\, dr\, du.$$

Using symmetry, the right-hand side is less than or equal to

$$\int_0^s \int_0^t \frac{p_1(0)}{(r+u)^{d/\beta}}\, dr\, du$$

with equality when $x_0 = y_0$. Some routine calculus completes the proof. □

**3. Renormalized self-intersection local times.** Let $X_t$ be a symmetric stable process of order $\beta$ in $R^d$. For any random variable $Y$, we set $\{Y\}_0 = Y - E(Y)$. For each bounded Borel set $B \subseteq R^2_+$, let

$$\gamma_\varepsilon(B) = \left\{ \int_B \int f_\varepsilon(X_s - X_r)\, dr\, ds \right\}_0. \tag{3.1}$$

We set $\gamma_{t,\varepsilon} = \gamma_\varepsilon(B_t)$, where $B_t = \{(r,s) \in R^2_+ | 0 \le r \le s \le t\}$.

Using the scaling $X_{\lambda s} \stackrel{d}{=} \lambda^{1/\beta} X_s$ and $f_{\lambda\varepsilon}(x) = \frac{1}{\lambda^d} f_\varepsilon(x/\lambda)$, we have

$$\gamma_\varepsilon(B) \stackrel{d}{=} \lambda^{-(2-d/\beta)} \gamma_{\lambda^{1/\beta}\varepsilon}(\lambda B). \tag{3.2}$$

THEOREM 7. *Let $X_t$ be a symmetric stable process of order $\beta > 2d/3$ in $R^d$. Then for all $\rho > 0$ sufficiently small, we can find some $\theta > 0$ such that*

$$\sup_{\varepsilon,\varepsilon',t>0} E\left( \exp\left\{ \theta \left| \frac{\gamma_{t,\varepsilon} - \gamma_{t,\varepsilon'}}{|\varepsilon - \varepsilon'|^\rho t^{2-(d+\rho)/\beta}} \right|^{\beta/(d+\rho)} \right\} \right) < \infty. \tag{3.3}$$

PROOF. Taking $\lambda = 1/t$ and $B = B_t$ in (3.2), we see that it suffices to prove (3.3) when $t = 1$. We adapt a technique pioneered by Le Gall [17].

Let

$$A^n_k = [(2k-2)2^{-n}, (2k-1)2^{-n}] \times [(2k-1)2^{-n}, (2k)2^{-n}]. \tag{3.4}$$

Note that $B_1 = \bigcup_{n=1}^\infty \bigcup_{k=1}^{2^{n-1}} A^n_k$ so that, for any $\varepsilon > 0$,

$$\gamma_{1,\varepsilon} = \sum_{n=1}^\infty \sum_{k=1}^{2^{n-1}} \gamma_\varepsilon(A^n_k). \tag{3.5}$$

We will use the following lemma whose proof is given at the end of this section.

LEMMA 1. *Let $0 < p \le 1$ and let $\{Y_k(\zeta)\}_{k\ge1}$ be a family (indexed by $\zeta$) of sequences of i.i.d. real valued random functions such that $E(Y_k(\zeta)) = 0$ and*

$$\lim_{\theta\to0} \sup_\zeta E e^{\theta|Y_1(\zeta)|^p} = 1. \tag{3.6}$$

*Then for some $\lambda > 0$,*

$$\sup_{n,\zeta} E \exp\left\{ \lambda \left| \sum_{k=1}^n Y_k(\zeta)/\sqrt{n} \right|^p \right\} < \infty. \tag{3.7}$$

By (2.4), for some $\rho > 0$,

$$\lim_{\theta\to 0}\sup_{\varepsilon,\varepsilon'>0} E\left(\exp\left\{\theta\left|\frac{\gamma_\varepsilon(A_1^1)-\gamma_{\varepsilon'}(A_1^1)}{|\varepsilon-\varepsilon'|^\rho}\right|^{\beta/(d+\rho)}\right\}\right)=1. \tag{3.8}$$

Hence, by Lemma 1, for some $\lambda > 0$,

$$\begin{aligned} e^\phi := \sup_{N,\varepsilon,\varepsilon'>0}\Bigg(E\Bigg(\exp\Bigg\{\lambda\Bigg|\sum_{k=1}^{2^{N-1}}\{\gamma_\varepsilon(2^{(N-1)}A_k^N)-\gamma_{\varepsilon'}(2^{(N-1)}A_k^N)\}\\ \times(2^{(N-1)/2}|\varepsilon-\varepsilon'|^\rho)^{-1}\Bigg|^{\beta/(d+\rho)}\Bigg\}\Bigg)\Bigg) \end{aligned} \tag{3.9}$$

is finite.

Since $\beta > \frac{2}{3}d$, for $\rho > 0$ sufficiently small,

$$a := \tfrac{3}{2}\beta/(d+\rho)-1>0. \tag{3.10}$$

Write

$$b_1=\lambda 2^{-a}\quad\text{and}\quad b_N=\lambda 2^{-a}\prod_{j=2}^{N}(1-2^{-aj}),\qquad N=2,3,\ldots. \tag{3.11}$$

Then for any integer $N \geq 1$, by Hölder's inequality,

$$\begin{aligned} \Psi_{\varepsilon,\varepsilon',N} &:= E\left(\exp\left\{b_N\left|\frac{\sum_{n=1}^{N}\sum_{k=1}^{2^{n-1}}\{\gamma_\varepsilon(A_k^n)-\gamma_{\varepsilon'}(A_k^n)\}}{|\varepsilon-\varepsilon'|^\rho}\right|^{\beta/(d+\rho)}\right\}\right)\\ &\le \Bigg(E\Bigg(\exp\Bigg\{\frac{b_N}{(1-2^{-aN})}\\ &\qquad\times\left|\frac{\sum_{n=1}^{N-1}\sum_{k=1}^{2^{n-1}}\{\gamma_\varepsilon(A_k^n)-\gamma_{\varepsilon'}(A_k^n)\}}{|\varepsilon-\varepsilon'|^\rho}\right|^{\beta/(d+\rho)}\Bigg\}\Bigg)\Bigg)^{1-2^{-aN}}\\ &\quad\times\left(E\left(\exp\left\{b_N 2^{aN}\left|\frac{\sum_{k=1}^{2^{N-1}}\{\gamma_\varepsilon(A_k^N)-\gamma_{\varepsilon'}(A_k^N)\}}{|\varepsilon-\varepsilon'|^\rho}\right|^{\beta/(d+\rho)}\right\}\right)\right)^{2^{-aN}}. \end{aligned} \tag{3.12}$$

Taking $\lambda = 2^{N-1}$ in (3.2), we see that

$$\begin{aligned} &\sum_{k=1}^{2^{N-1}}\{\gamma_\varepsilon(A_k^N)-\gamma_{\varepsilon'}(A_k^N)\}\\ &\quad\stackrel{d}{=}2^{-(2-d/\beta)(N-1)}\\ &\qquad\times\sum_{k=1}^{2^{N-1}}\{\gamma_{\varepsilon 2^{(N-1)/\beta}}(2^{(N-1)}A_k^N)-\gamma_{2^{(N-1)/\beta}\varepsilon'}(2^{(N-1)}A_k^N)\}. \end{aligned} \tag{3.13}$$

Using (3.10), we note that

$$\left(2-\frac{d}{\beta}\right)-\frac{\rho}{\beta}-a\frac{(d+\rho)}{\beta}=\frac{1}{2}. \tag{3.14}$$

Hence,

$$
\begin{aligned}
&2^{aN}\left|\frac{\sum_{k=1}^{2^{N-1}}\{\gamma_\varepsilon(A_k^N)-\gamma_{\varepsilon'}(A_k^N)\}}{|\varepsilon-\varepsilon'|^\rho}\right|^{\beta/(d+\rho)}\\
&\quad\le 2^a\left|\frac{\sum_{k=1}^{2^{N-1}}\{\gamma_{\varepsilon 2^{(N-1)/\beta}}(2^{(N-1)}A_k^N)-\gamma_{\varepsilon' 2^{(N-1)/\beta}}(2^{(N-1)}A_k^N)\}}{2^{(N-1)/2}|\varepsilon 2^{(N-1)/\beta}-\varepsilon' 2^{(N-1)/\beta}|^\rho}\right|^{\beta/(d+\rho)}
\end{aligned}
\tag{3.15}
$$

in law. Using this, the finiteness of (3.9) and the fact that $b_N 2^a \le \lambda$ for the last line of (3.12), and (3.11) and the fact that $1-2^{-aN}<1$ for the second line of (3.12), we have that

$$
\Psi_{\varepsilon,\varepsilon',N}\le \Psi_{\varepsilon,\varepsilon',N-1}\exp\{\phi 2^{-aN}\}. \tag{3.16}
$$

Inductively,

$$
\Psi_{\varepsilon,\varepsilon',N}\le \exp\{\phi 2^{-a}(1-2^{-a})^{-1}\}.
$$

Letting $N\to\infty$, Theorem 7 follows by (3.5) and Fatou's lemma. □

It follows from Theorem 7 and Kolmogorov's continuity theorem that

$$
\gamma_t := \lim_{\varepsilon\to 0}\gamma_{\varepsilon,t} \tag{3.17}
$$

exists a.s. and in all $L^p$ spaces.

Furthermore, it follows from Theorem 7 that for some $\rho,\theta>0$,

$$
\sup_{\varepsilon,t>0} E\left(\exp\left\{\theta\left|\frac{\gamma_t-\gamma_{t,\varepsilon}}{\varepsilon^\rho t^{2-(d+\rho)/\beta}}\right|^{\beta/(d+\rho)}\right\}\right)<\infty. \tag{3.18}
$$

Note that, since for $\rho>0$ sufficiently small $\beta/(d+\rho)>1/2$, it follows that for any $\lambda,\delta>0$,

$$
\begin{aligned}
&E(\exp\{\lambda|\gamma_t-\gamma_{t,\varepsilon}|^{1/2}\})\\
&\quad\le e^{\lambda\delta t}+E(\exp\{\lambda|\gamma_t-\gamma_{t,\varepsilon}|^{1/2}\}\mathbb{1}_{\{|\gamma_t-\gamma_{t,\varepsilon}|\ge(\delta t)^2\}})\\
&\quad\le e^{\lambda\delta t}+E\left(\exp\left\{\lambda\left|\frac{\gamma_t-\gamma_{t,\varepsilon}}{(\delta t)^{2-(d+\rho)/\beta}}\right|^{\beta/(d+\rho)}\right\}\right).
\end{aligned}
\tag{3.19}
$$

Using (3.18), we conclude that, for any $\lambda>0$,

$$
\limsup_{\varepsilon\to 0}\limsup_{t\to\infty}\frac{1}{t}\log E(\exp\{\lambda|\gamma_t-\gamma_{t,\varepsilon}|^{1/2}\})=0. \tag{3.20}
$$

For later reference we note that arguments similar to those used in proving Theorem 7 show that, for some $\theta>0$,

$$
\sup_{t>0} E\left(\exp\left\{\theta\left|\frac{\gamma_t}{t^{2-d/\beta}}\right|^{\beta/d}\right\}\right)<\infty. \tag{3.21}
$$

(In fact, by scaling, we only need this for $t=1$.)

PROOF OF LEMMA 1. Let $\psi_p(x) = e^{x^p} - 1$ for large $x$ and linear near the origin so that $\psi_p(x)$ is convex. We use $\|\cdot\|_{\psi_p}$ to denote the norm of the Orlicz space $L_{\psi_p}$ with Young's function $\psi_p$. Assumption (3.6) implies that, for some $M < \infty$,

$$\sup_{\zeta} \|Y_1(\zeta)\|_{\psi_p} \le M. \tag{3.22}$$

By Theorem 6.21 of [13], if $\xi_k$ are i.i.d. copies of a mean zero random variable $\xi_1 \in L_{\psi_p}$, then for some constant $K_p$, depending only on $p$,

$$\left\| \sum_{k=1}^{n} \xi_k \right\|_{\psi_p} \le K_p \left( \left\| \sum_{k=1}^{n} \xi_k \right\|_{L_1} + \left\| \max_{1\le k\le n} |\xi_k| \right\|_{\psi_p} \right).$$

Using Proposition 4.3.1 of [11], for some constant $C_p$, depending only on $p$,

$$\left\| \max_{1\le k\le n} |\xi_k| \right\|_{\psi_p} \le C_p (\log n) \|\xi_1\|_{\psi_p}.$$

Since the $\xi_k$ are i.i.d. and mean zero,

$$\left\| \sum_{k=1}^{n} \xi_k \right\|_{L_1} \le \left\| \sum_{k=1}^{n} \xi_k \right\|_{L_2} \le \sqrt{n} \|\xi_1\|_{L_2}.$$

Thus, we have

$$\left\| \sum_{k=1}^{n} \xi_k / \sqrt{n} \right\|_{\psi_p} \le D_p \left( \|\xi_1\|_{L_2} + \frac{\log n}{\sqrt{n}} \|\xi_1\|_{\psi_p} \right)$$

for some constant $D_p$, depending only on $p$. Lemma 1 follows immediately from this. $\square$

## 4. Large deviations for renormalized self-intersection local times.

Let

$$\mathcal{E}_\psi(f,f) := \int_{R^d} \psi(p) |\widehat{f}(p)|^2 \, dp \tag{4.1}$$

and set

$$\mathcal{F}_\psi = \{ f \in L^2(R^d) \mid \|f\|_2 = 1, \mathcal{E}_\psi(f,f) < \infty \}. \tag{4.2}$$

The following lemma is proven is Section 2 of [8].

LEMMA 2. *If $\beta > d/2$, then for any $\lambda > 0$,*

$$M_\psi(\lambda) := \sup_{f \in \mathcal{F}_\psi} \{ \lambda \|f\|_4^2 - \mathcal{E}_\psi(f,f) \} < \infty \tag{4.3}$$

*and*

$$M_\psi(\lambda) = \lambda^{2\beta/(2\beta - d)} M_\psi(1). \tag{4.4}$$

*Furthermore,*

$$\kappa_\psi := \inf\{C \mid \|f\|_{2p} \le C\|f\|_2^{1-d/2\beta}[\mathcal{E}_\psi^{1/2}(f,f)]^{d/2\beta}\} < \infty \tag{4.5}$$

*and*

$$M_\psi(1) = \frac{2\beta - d}{d}\left(\frac{d\kappa_\psi^2}{2\beta}\right)^{2\beta/(2\beta-d)}. \tag{4.6}$$

We write $M_\psi = M_\psi(1)$ and let

$$K_\psi = \frac{d}{\beta}\left(\frac{2\beta-d}{2\beta M_\psi}\right)^{(2\beta-d)/d}. \tag{4.7}$$

PROOF OF THEOREM 1. We show that if $X_t$ is a symmetric stable process of order $\beta > 2d/3$ in $R^d$, then

$$\lim_{t\to\infty}\frac{1}{t}\log P(\gamma_t \ge t^2) = -2^{\beta/d-1}K_\psi. \tag{4.8}$$

[This defines $a_\psi$ of (1.7).]

Let $h$ be a positive, symmetric function in the Schwarz class $\mathcal{S}(R^d)$ with $\int h\,dx = 1$, and note that $f = h * h$ has the same properties and $f_\varepsilon = h_\varepsilon * h_\varepsilon$. Using this, observe that

$$\begin{aligned}
&\int_0^t\int_0^s f_\varepsilon(X_s - X_r)\,dr\,ds\\
&\qquad = \tfrac{1}{2}\int_0^t\int_0^t f_\varepsilon(X_s - X_r)\,dr\,ds\\
&\qquad = \tfrac{1}{2}\int_{R^d}\left(\int_0^t h_\varepsilon(X_s - x)\,ds\right)^2 dx,
\end{aligned} \tag{4.9}$$

hence, by Theorem 5 of [8], for any $\lambda > 0$,

$$\begin{aligned}
&\lim_{t\to\infty}\frac{1}{t}\log E\exp\left\{\lambda\left(\int_0^t\int_0^s f_\varepsilon(X_s - X_r)\,dr\,ds\right)^{1/2}\right\}\\
&\qquad = \lim_{t\to\infty}\frac{1}{t}\log E\exp\left\{\frac{\lambda}{\sqrt{2}}\left(\int_{R^d}\left(\int_0^t h_\varepsilon(X_s - x)\,ds\right)^2 dx\right)^{1/2}\right\}\\
&\qquad = \sup_{g\in\mathcal{F}_\psi}\left\{\frac{\lambda}{\sqrt{2}}\left(\int_{R^d}|(g^2 * h_\varepsilon)(x)|^2\,dx\right)^{1/2} - \mathcal{E}_\psi(g,g)\right\}.
\end{aligned} \tag{4.10}$$

For each fixed $\varepsilon > 0$,

$$\begin{aligned}
&E\left(\int_0^t\int_0^s f_\varepsilon(X_s - X_r)\,dr\,ds\right)\\
&\qquad = \int_{R^d}\int_0^t\int_0^s E(e^{ip\cdot(X_s - X_r)})\,dr\,ds\,\widehat{f}(\varepsilon p)\,dp\\
&\qquad = \int_{R^d}\int_0^t\int_0^s e^{-(s-r)\psi(p)}\,dr\,ds\,\widehat{f}(\varepsilon p)\,dp\\
&\qquad \le Ct\int_{R^d}\frac{1}{|p|^\beta}\widehat{f}(\varepsilon p)\,dp = O(t)
\end{aligned} \tag{4.11}$$

if $\beta < d$. [When $\beta = d$, we can easily obtain $O(t^{1+\delta})$ for any $\delta > 0$.] Using (3.20), we conclude that for any $\lambda > 0$,

$$\limsup_{\varepsilon\to 0}\limsup_{t\to\infty}\frac{1}{t}\log E\biggl(\exp\biggl\{\lambda\biggl|\gamma_t-\int_0^t\int_0^s f_\varepsilon(X_s-X_r)\,dr\,ds\biggr|^{1/2}\biggr\}\biggr)=0. \tag{4.12}$$

Hence, using (4.10) together with the argument used to take the $\varepsilon \to 0$ limit in [8] and then recalling (4.4),

$$\begin{aligned}
&\lim_{t\to\infty}\frac{1}{t}\log E\exp\{\lambda|\gamma_t|^{1/2}\}\\
&\quad=\lim_{\varepsilon\to 0}\sup_{g\in\mathcal{F}_\psi}\biggl\{\frac{\lambda}{\sqrt{2}}\biggl(\int_{R^d}|(g^2*h_\varepsilon)(x)|^2\,dx\biggr)^{1/2}-\mathcal{E}_\psi(g,g)\biggr\}\\
&\quad=\sup_{g\in\mathcal{F}_\psi}\biggl\{\frac{\lambda}{\sqrt{2}}\biggl(\int_{R^d}g^4(x)\,dx\biggr)^{1/2}-\mathcal{E}_\psi(g,g)\biggr\}\\
&\quad=\biggl(\frac{\lambda}{\sqrt{2}}\biggr)^{2\beta/(2\beta-d)}M_\psi.
\end{aligned} \tag{4.13}$$

By the Gärtner–Ellis theorem ([9], Theorem 2.3.6)

$$\begin{aligned}
&\lim_{t\to\infty}\frac{1}{t}\log P(|\gamma_t|\ge t^2)\\
&\quad=-\sup_{\lambda>0}\biggl\{\lambda-\biggl(\frac{\lambda}{\sqrt{2}}\biggr)^{2\beta/(2\beta-d)}M_\psi\biggr\}\\
&\quad=-2^{\beta/d-1}\frac{d}{\beta}\biggl(\frac{2\beta-d}{2\beta M_\psi}\biggr)^{(2\beta-d)/d}.
\end{aligned} \tag{4.14}$$

On the other hand, writing $\gamma_t = \gamma_t^+ - \gamma_t^-$ and using the positivity of $\int_0^t\int_0^s f_\varepsilon(X_s - X_r)\,dr\,ds$ and (4.12), we have that for any $\lambda$,

$$\limsup_{t\to\infty}\frac{1}{t}\log E(\exp\{\lambda|\gamma_t^-|^{1/2}\})=0. \tag{4.15}$$

Theorem 1 then follows. $\square$

## 5. The lower tail; $\beta < d$.

Proof of Theorem 2 when $\beta < d$. For each bounded Borel set $A \subseteq R_+^2$, we set $\gamma(A) = \lim_{\varepsilon\to 0}\gamma_\varepsilon(A)$, recall (3.1). This limit is known to exist. Let $\Gamma([s,t]) := \gamma(\{(u,v)|s \le u \le v \le t\})$ and with $[0,s;s,t] = \{(u,v)|0 \le u \le s \le v \le t\}$ note that $\gamma([0,s;s,t]) \overset{d}{=} \{\alpha_{s,t-s}\}_0$. Thus, for any positive $s$ and $t$,

$$\begin{aligned}
\gamma_{s+t}&=\gamma_s+\Gamma([s,s+t])+\gamma([0,s];[s,s+t])\\
&\ge\gamma_s+\Gamma([s,s+t])-E\alpha([0,s];[s,s+t]).
\end{aligned} \tag{5.1}$$

Note that $\gamma_s \in \mathcal{F}_s = \sigma(X_r, 0 \le r \le s)$, $\Gamma([s, s+t])$ is independent of $\mathcal{F}_s$, and $\Gamma([s, s+t])$ has the same distribution as $\gamma_t$. Define

$$Z_t = c_\psi t^{2-d/\beta} - \gamma_t, \qquad Z_{s,t} = c_\psi t^{2-d/\beta} - \Gamma([s, s+t]). \tag{5.2}$$

By the above, $\{Z_{s,t}; t \ge 0\}$ is independent of $\{Z_u; u \le s\}$ and we have $\{Z_{s,t}; t \ge 0\} \stackrel{d}{=} \{Z_t; t \ge 0\}$. Using (5.1) and Theorem 6, we have that for any $s, t > 0$,

$$Z_{s+t} \le Z_s + Z_{s,t}. \tag{5.3}$$

Given $a > 0$, define

$$\tau_a = \inf\{s; Z_s \ge a\}.$$

By continuity, $Z_{\tau_a} = a$ on $\tau_a < \infty$. Let

$$\phi(h) = \sup_{\substack{0 \le s,t \le 1 \\ |t-s| \le h}} |Z_t - Z_s|. \tag{5.4}$$

Fix $a, b, n > 0$ and $0 < \delta < a, b$,

$$\begin{aligned}
&P\Big(\sup_{t \le 1} Z_t \ge a + b,\, \phi(1/n) \le \delta\Big) \\
&\quad = \sum_{j=0}^{n-2} P\Big(\sup_{t \le 1} Z_t \ge a + b, \phi(1/n) \le \delta, j/n \le \tau_a < (j+1)/n\Big) \\
&\quad \le \sum_{j=0}^{n-2} P\Big(\sup_{t \le 1} Z_{(j+1)/n,t} \ge b - \delta, j/n \le \tau_a < (j+1)/n\Big) \\
&\quad = \sum_{j=0}^{n-2} P\Big(\sup_{t \le 1} Z_{(j+1)/n,t} \ge b - \delta\Big) P(j/n \le \tau_a < (j+1)/n) \\
&\quad \le P\Big(\sup_{t \le 1} Z_t \ge a\Big) P\Big(\sup_{t \le 1} Z_t \ge b - \delta\Big).
\end{aligned} \tag{5.5}$$

Using the continuity of $Z_s$ and first taking $n \to \infty$ and then $\delta \to 0$, we obtain

$$P\Big(\sup_{t \le 1} Z_t \ge a + b\Big) \le P\Big(\sup_{t \le 1} Z_t \ge a\Big) P\Big(\sup_{t \le 1} Z_t \ge b\Big). \tag{5.6}$$

Hence, there is $c > 0$ such that for some $\lambda_0 < \infty$,

$$P\Big(\sup_{t \le 1} Z_t \ge \lambda\Big) \le e^{-c\lambda} \qquad \forall\, \lambda > \lambda_0, \tag{5.7}$$

so that

$$E \exp\Big\{c_0 \sup_{t \le 1} Z_t\Big\} < \infty \tag{5.8}$$

for some $c_0 > 0$. Then by the sub-additivity (5.3) and what we have just proven, there is $c_0 > 0$ such that

$$E\exp\left\{c_0 \sup_{t\le n} Z_t\right\} \le \left(E\exp\left\{c_0 \sup_{t\le 1} Z_t\right\}\right)^n < \infty$$

for all $n$. Then by the scaling (1.10), we see that (5.8) holds for all $c_0 > 0$. Therefore, we have

$$E\exp\left\{c\sup_{t\le n}\{-\gamma_t\}\right\} < \infty \qquad \forall\, c, n > 0. \tag{5.9}$$

Setting now

$$a_\lambda(t) = \log(E\exp\{\lambda Z_t\}),$$

by the sub-additivity (5.3), we have that for any positive $s, t, \lambda$,

$$a_\lambda(s+t) \le a_\lambda(s) + a_\lambda(t). \tag{5.10}$$

Consequently,

$$\lim_{t\to\infty} \frac{1}{t} a_\lambda(t) = \inf_{t\ge 1}\left\{\frac{1}{t}a_\lambda(t)\right\} := L_\lambda < \infty, \tag{5.11}$$

where the last inequality follows from (5.9). Note that

$$a_\lambda(t) = \lambda c_\psi t^{2-d/\beta} + \log(E\exp\{-\lambda\gamma_t\}),$$

with $2 - d/\beta < 1$, so that (5.11) implies that for any $\lambda > 0$,

$$\lim_{t\to\infty} \frac{1}{t}\log(E\exp\{-\lambda\gamma_t\}) = L_\lambda < \infty. \tag{5.12}$$

It follows from Theorem 8, immediately following, that $L_{\lambda_0} > 0$ for some $0 < \lambda_0 < \infty$. Using the scaling (1.10), it follows from (5.12) that for any $\lambda > 0$,

$$\lim_{t\to\infty} \frac{1}{t}\log(E\exp\{-\lambda\gamma_t\}) = \lambda^{\beta/(2\beta-d)}\lambda_0^{-\beta/(2\beta-d)}L_{\lambda_0}. \tag{5.13}$$

It then follows by the Gärtner–Ellis theorem, compare (4.13) and (4.14), that

$$\lim_{t\to\infty} t^{-1}\log P(-\gamma_t \ge t) = -b_\psi, \tag{5.14}$$

with

$$b_\psi = \left(\frac{d-\beta}{\beta}\right)\left(\frac{2\beta-d}{\beta L_{\lambda_0}}\right)^{(2\beta-d)/(d-\beta)} \lambda_0^{\beta/(d-\beta)}.$$

Note that it follows from (5.13) that $\lambda_0^{-\beta/(2\beta-d)}L_{\lambda_0}$ is independent of the particular $\lambda_0$ chosen so the same will be true of $b_\psi$. This will complete the proof of Theorem 2 when $\beta < d$. □

THEOREM 8. *Let $X_t$ be a symmetric stable process of order $\beta \in (2d/3, d)$ in $R^d$. There exist constants $c_1, c_2 > 0$ such that*

$$P(-\gamma_n \geq c_1 n) \geq c_2^n. \tag{5.15}$$

The idea of the proof is the following. Let $\varepsilon$ be small, $M = \varepsilon^{-1}$ and $Q_k$ the square with one diagonal going from the point $(Mk - 4\varepsilon, 0)$ to the point $(M(k+1) + 4\varepsilon, 0)$. By scaling and some easy estimates, we show that, for each $k$, there is probability on the order of $\varepsilon$ to a power that $X_t$ lies in $Q_k$ when $t \in [k, k+1]$ and also the renormalized self-intersection local time of that portion of the path of $X$ is not too small. Provided the intersection local times between consecutive portions of the path are not too large, we can then use the Markov property $n$ times to obtain the result of Theorem 8. The intersection local time of consecutive portions of the path may be viewed as the intersection local time of two independent stable processes. We use the representation of this intersection local time as an additive functional along the lines of [3] to obtain a suitable upper bound on its size, except for a set whose probability decreases faster than any power of $\varepsilon$. We then take $\varepsilon$ sufficiently small, but fixed.

PROOF OF THEOREM 8. Let $A(I; J)$ denote the intersection local time between $X(I)$ and $X(J)$, where $X(I) = \{X_s : s \in I\}$ for an interval $I$ and let $\Gamma(I)$ denote the renormalized self-intersection local time of $X(I)$. $\varepsilon < 1/4$ will be chosen later. Set $M = \varepsilon^{-1}$. First of all, $-\Gamma([0,1])$ has mean 0 and is not identically zero. So there exist positive constants $\kappa_1, \kappa_2$ not depending on $\varepsilon$ such that

$$P(-\Gamma([0,1]) > \kappa_1) > \kappa_2.$$

By scaling,

$$P(-\Gamma([\varepsilon^2, 1-\varepsilon^2]) > \kappa_1/2) > \kappa_2.$$

If we choose $\varepsilon$ small enough, by the fact that the paths of $X_t$ are right continuous with left limits,

$$P\Big(\sup_{\varepsilon^2 \leq s \leq 1-\varepsilon^2} |X_s - X_{\varepsilon^2}| > M/2\Big) \leq \kappa_2/2.$$

Therefore, if

$$E_1 = \Big\{-\Gamma([\varepsilon^2, 1-\varepsilon^2]) > \kappa_1/2, \sup_{\varepsilon^2 \leq s \leq 1-\varepsilon^2} |X_s - X_{\varepsilon^2}| \leq M/2\Big\},$$

then

$$P(E_1) \geq \kappa_2/2.$$

Let $B(x,r)$ denote the open ball in $R^d$ of radius $r$ centered at $x$. Let $S_k = B((Mk,0),\varepsilon^2)$, that is, the ball with center at the point $(Mk,0)$ and radius $\varepsilon$, and let $Q_k$ be the square which has one diagonal going from $(Mk-4\varepsilon,0)$ to $(M(k+1)+4\varepsilon,0)$. Let $z_k$ be the center of $Q_k$, that is, $z_k = (M(k+\frac{1}{2}),0)$. Let

$$E_2 = \{X_{\varepsilon^2} \in B(z_k,1) \text{ and } X_s \in Q_k \text{ for } s \in [0,\varepsilon^2]\}.$$

Let

$$E_3 = \{X_{\varepsilon^2} \in S_{k+1} \text{ and } X_s \in Q_k \text{ for } s \in [0,\varepsilon^2]\}.$$

As usual, we use $P^x$ for the probability when our process $X$ is started at $x$.

LEMMA 3. (a) *There exists $c_3$ such that if $x \in S_k$ and $\varepsilon$ is sufficiently small, then*

$$P^x(E_2) \geq c_3\varepsilon^{4+\beta}.$$

(b) *If $x \in B(z_k, M/2)$ and $\varepsilon$ is sufficiently small, then*

$$P^x(E_3) \geq c_3\varepsilon^{6+\beta}.$$

PROOF. (a) Let $\tau = \inf\{t : |X_t - X_0| > \varepsilon/2\}$. By scaling and the fact that $\beta > 1$, we have $P(\sup_{s\leq\varepsilon^2}|X_s - X_0| > \varepsilon/2) \to 0$ as $\varepsilon \to 0$. So by taking $\varepsilon$ small enough, we may assume that

$$P^x(\tau \leq \varepsilon^2) \leq 1/2$$

for all $x$.

By the Lévy system formula for right continuous stable processes (see [4], Proposition 2.3, e.g.),

$$\begin{aligned} &P^x(X_{\tau\wedge\varepsilon^2} \in B(z_k,1/2)) \\ &\quad \geq E^x \sum_{s\leq\tau\wedge\varepsilon^2} \mathbb{1}_{(X_{s-}\in B((Mk,0),\varepsilon/2))}\mathbb{1}_{(X_s\in B(z_k,1/2))} \\ &\quad = E^x \int_0^{\tau\wedge\varepsilon^2}\int_{B(z_k,1/2)} n(X_s,z)\,dz\,ds, \end{aligned} \tag{5.16}$$

where $n(y,z) = c_4|y-z|^{-2-\beta}$. Since $n(y,z)$ is bounded below by $c_4M^{-2-\beta}$ if $y \in B((Mk,0),2\varepsilon)$ and $z \in B(z_k,1/2)$, we see

$$\begin{aligned} &P^x(X_{\tau\wedge\varepsilon^2} \in B(z_k,1/2)) \\ &\quad \geq c_4\varepsilon^{2+\beta}E^x[\tau\wedge\varepsilon^2] \geq c_4\varepsilon^{2+\beta}E^x[\varepsilon^2;\tau > \varepsilon^2] \\ &\quad = c_4\varepsilon^{2+\beta}\varepsilon^2P^x(\tau > \varepsilon^2) \geq c_3\varepsilon^{4+\beta}/2. \end{aligned} \tag{5.17}$$

We noted in the first paragraph of the proof that there is probability at least $1/2$ that $X_t$ moves no more than $\varepsilon/2$ in time $\varepsilon^2$. So by using the strong Markov property at time $\tau$, there is probability at least $c_4\varepsilon^{4+\beta}/4$ that $X_t$ exits $S_k$ by time $\varepsilon^2$, jumps to $B(z_k,1/2)$, and then stays in $B(z_k,1)$ until time $\tau+\varepsilon^2$. But this event is contained in $E_2$.

(b) The proof of (b) is similar. Using the Lévy system formula,

$$\begin{aligned} P^x&(X_{\tau\wedge\varepsilon^2}\in B(M((k+1),0),\varepsilon/2)) \\ &\geq E^x\int_0^{\tau\wedge\varepsilon^2}\int_{B((M(k+1),0),\varepsilon/2)} n(X_s,z)\,dz\,ds. \end{aligned}$$

This, in turn, is greater than or equal to

$$c_5\varepsilon^2M^{-2-\beta}E^x[\tau\wedge\varepsilon^2]\geq c_6\varepsilon^{6+\beta}.$$

We chose $\varepsilon$ so that the probability that $X_t$ moves no more than $\varepsilon/2$ in time $\varepsilon^2$ is at least $1/2$. Using the strong Markov property at time $\tau$, there is probability at least $c_6\varepsilon^{6+\beta}/2$ that the process exits $B(x,\varepsilon/2)$ by time $\varepsilon^2$, jumps to $B((M(k+1),0),\varepsilon/2)$, and then moves no more than $\varepsilon/2$ in time $\varepsilon^2$. This event is contained in $E_3$, and (b) follows. This completes the proof of Lemma 3. $\square$

Let

$$E_3'=E_3\circ\theta_{1-\varepsilon^2}=\{X_1\in S_{k+1}\text{ and }X_s\in Q_k\text{ for }s\in[1-\varepsilon^2,1]\}.$$

Using Lemma 3 and the Markov property at times $\varepsilon^2$ and $1-\varepsilon^2$,

$$P^x(E_1\cap E_2\cap E_3')\geq c_3^2\varepsilon^{10+2\beta}\kappa_2/2. \tag{5.18}$$

Let

$$\begin{aligned} E_4&=\{\Gamma[0,\varepsilon^2]>\kappa_1/16\},\\ E_5&=\{\Gamma[1-\varepsilon^2,1]>\kappa_1/16\},\\ E_6&=\{A([0,\varepsilon^2];[\varepsilon^2,1])>\kappa_1/16\},\\ E_7&=\{A([0,1-\varepsilon^2];[1-\varepsilon^2,1])>\kappa_1/16\}. \end{aligned} \tag{5.19}$$

LEMMA 4. *There exist $c_7,c_8$ and $b$ not depending on $\varepsilon$ such that*

$$P(E_4)+P(E_5)+P(E_6)+P(E_7)\leq c_7e^{-c_8/\varepsilon^b}.$$

PROOF. The estimates for $E_4$ and $E_5$ follow from the scaling (1.10) and (1.14). By (2.16),

$$P(A([0,1];[1,1+a])>\lambda)\leq c_9e^{-c_{10}\lambda^{\beta/d}/a^{\beta/d-1/2}}. \tag{5.20}$$

This and scaling give us the desired estimates for $E_6$ and $E_7$. This completes the proof of Lemma 4. $\square$

Recall that the *occupation measure* $\mu_T^X$ is defined as

$$\mu_t^X(A) = \int_0^t \mathbb{1}_A(X_s)\, ds$$

for all Borel sets $A \subseteq R^d$. If $p_s(x)$ is the probability density function for $X_s$ and $u(x) = \int_0^\infty p_s(x)\, ds$ is the 0-potential density for $X$, it is easily checked that

$$E^x(\{\mu_\infty^X(A)\}^n) = n! \int \prod_{j=1}^n u(x_i - x_{i-1}) \mathbb{1}_A(x_i)\, dx_i, \tag{5.21}$$

where $x_0 = x$. Hence, if

$$c_A = \sup_x \int u(x-y) \mathbb{1}_A(y)\, dy, \tag{5.22}$$

we have that $\sup_x E^x(\{\mu_\infty^X(A)\}^n) \le n! c_A^n$ and, thus,

$$\sup_x E^x(\exp\{\mu_\infty^X(A)/2c_A\}) \le 2$$

so that, by Chebyshev,

$$\sup_x P^x(\mu_\infty^X(A) \ge 2\lambda c_A) \le 2e^{-\lambda}. \tag{5.23}$$

LEMMA 5. *Let $\delta \in (0, 2\beta - 2)$ and $M > 2$. There exist constants $c_{11}$ and $c_{12}$ depending only on $M$ and $\delta$ such that*

$$P\left(\sup_{|x| \le M, 0 < r \le 1} \frac{\mu_\infty^X(B(x,r))}{r^{\beta-\delta}} > \lambda\right) \le c_{11} M^2 e^{-c_{12}\lambda}. \tag{5.24}$$

PROOF. First fix $x$ and $r$. Since $u(y-z) \le c_{13}|y-z|^{\beta-2}$, using symmetry, $c_{B(x,r)}$ is bounded by

$$\int_{B(x,r)} c_{13}|x-z|^{\beta-2}\, dz = c_{14} r^\beta.$$

Applying (5.23),

$$P(\mu_\infty^X(B(x,r)) > \lambda r^{\beta-\delta}) \le 2e^{-c_{15}\lambda r^{-\delta}}. \tag{5.25}$$

Suppose now that $\mu_\infty^X(B(x,r)) > \lambda r^{\beta-\delta}$ for some $|x| \le M$ and some $r \in (0,1)$. Choose $k$ such that $2^{-k-1} \le r < 2^{-k}$ and choose $x'$ so that both coordinates of $x'$ are integer multiples of $2^{-k}$ and $|x - x'| \le 2^{-k+1}$. Therefore,

$$\mu_\infty^X(B(x', 2^{-k+3})) > c_{16}\lambda(2^{-k+3})^{\beta-\delta},$$

where $c_{16}$ does not depend on $k$.

Since there are at most $c_{17}M^2 2^{2k}$ points in $B(0,2M)$ such that both coordinates are integer multiples of $2^{-k}$, then if $2^{-k-1} \le r < 2^{-k}$,

$$P\left(\sup_{|x|\le M} \frac{\mu_\infty^X(B(x,r))}{r^{\beta-\delta}} > c_{16}\lambda\right) \le c_{18} 2^{2k} M^2 e^{-c_{18}\lambda 2^{-\delta k}}. \tag{5.26}$$

Summing the right-hand side of (5.26) over $k$ from $-4$ to $\infty$ yields the right-hand side of (5.24). This completes the proof of Lemma 5. $\square$

By Lemma 5, it follows that

$$P\left(\sup_{|x|\le M, 0<r\le 1} \frac{\mu_\infty^X(B(x,r))}{r^{\beta-\delta}} > \kappa_1 \log^2(1/\varepsilon)/8\right) \le c_3^2 \varepsilon^{10+2\beta}\kappa_2/4 \tag{5.27}$$

if $\varepsilon$ is small enough.

Let $\mu_{t,t'}^X(A) = \int_t^{t'} \mathbb{1}_A(X_s)\,ds$, set

$$D_k = \Big\{X_k \in S_k, X_{k+1} \in S_{k+1}, \text{ and for } k \le s \le k+1, X_s \in Q_k,$$
$$-\Gamma[0,1] \ge \kappa_1/4, \sup_{|x|\le M, 0<r\le 1} \frac{\mu_{k,k+1}^X(B(x,r))}{r^{\beta-\delta}} \le \kappa_1 \log^2(1/\varepsilon)/8\Big\},$$

and recall that

$$\mathcal{F}_k = \sigma(X_v; v \le k).$$

By (5.18), Lemma 4, (5.27) and the Markov property,

$$P(D_k|\mathcal{F}_k) \ge c_{19}\varepsilon^{10+2\beta}\kappa_2/4 \qquad \text{on } D_{k-1}. \tag{5.28}$$

Let

$$F_k = \{A([k-1,k];[k,k+1]) \le \kappa_1/8\}, \qquad F_0 = \Omega,$$

and

$$L_k = D_k \cap F_k.$$

Lemma 6. *Let $\delta \in (0, 2\beta-2)$. We have*

$$P(F_k^c \cap D_k|\mathcal{F}_k) \le c_{20} e^{-c_{21}/\varepsilon^{2\beta-2-\delta}} \qquad \text{on } \bigcap_{j=1}^{k-1} L_j. \tag{5.29}$$

Proof. When $k=0$, there is nothing to prove, so let us suppose $k \ge 1$. As before, $A([k-1,k];[k,k+1])$ has the distribution of $\alpha_1$, and using the properties of $D_{k-1}, D_k$ and the Markov property, we have, recalling (2.1),

$$\begin{aligned} &P(F_k^c \cap D_k|\mathcal{F}_k) \\ &\quad \le \sup_{x\in S_k, X'\in D'_k} P_X^x\left(\lim_{\rho\to 0}\int_0^1\int_0^1 f_\rho(X_s - X'_r)\mathbb{1}_{Q_k}(X_s)\,dr\,ds \ge \kappa_1/8\right), \end{aligned} \tag{5.30}$$

where $P_X^x$ denotes probability with respect to the process $X$, while the independent process $X'$ is fixed, and

$$D_k' = \Big\{\mu_1^{X'}(\cdot) \text{ is supported on } Q_{k-1}, \\ \sup_{|x|\le M, 0<r\le 1} \frac{\mu_1^{X'}(B(x,r))}{r^{\beta-\delta}} \le \kappa_1 \log^2(1/\varepsilon)/8\Big\}.$$

In (5.30) we can and will take $f$ to be supported in $B(0,1)$. To bound the probability in (5.30), we note that

$$\lim_{\rho\to 0}\int_0^1\int_0^1 f_\rho(X_s - X_r')\mathbb{1}_{Q_k}(X_s)\,dr\,ds \\ \le \liminf_{\rho\to 0}\int_0^\infty\int_0^1 f_\rho(X_s - X_r')\mathbb{1}_{Q_k}(X_s)\,dr\,ds$$

and, by Fatou,

$$\begin{aligned} &E_X^x\bigg(\Big\{\liminf_{\rho\to 0}\int_0^\infty\int_0^1 f_\rho(X_s - X_r')\mathbb{1}_{Q_k}(X_s)\,dr\,ds\Big\}^n\bigg) \\ &\quad\le n!\liminf_{\rho\to 0}\int_{[0,1]^{nd}}\int_{R^{nd}}\prod_{j=1}^n u(x_i - x_{i-1}) f_\rho(x_i - X_{r_i}')\mathbb{1}_{Q_k}(x_i)\,dx_i\,dr_i \qquad (5.31)\\ &\quad= n!\liminf_{\rho\to 0}\int_{R^{nd}}\prod_{j=1}^n u(x_i - x_{i-1})\mathbb{1}_{Q_k}(x_i)\,d\mu_{1,\rho}^{X'}(x_i), \end{aligned}$$

with $x_0 = x$ and $d\mu_{1,\rho}^{X'}(x) = \int_0^1 f_\rho(x - X_r')\,dr\,dx$. As in the proof of (5.23), it then follows that $P(F_k^c \cap D_k|\mathcal{F}_k) \le c_{22}e^{-c_{23}/\bar{c}}$, where

$$\bar{c} = \sup_{0<\rho<\varepsilon}\ \sup_{x\in Q_{k-1}\cap Q_k, X'\in D_k'}\int_{R^d} u(y-x)\mathbb{1}_{Q_k}(y)\,d\mu_{1,\rho}^{X'}(y). \qquad (5.32)$$

It is easily checked that if $X' \in D_k'$, then uniformly in $\rho < \varepsilon$ and $0 < r \le 1 - \varepsilon$,

$$\sup_{|x|\le M-\varepsilon} \mu_{1,\rho}^{X'}(B(x,r)) \le cr^{\beta-\delta}\log^2(1/\varepsilon) \qquad (5.33)$$

and $\mu_{1,\rho}^{X'}$ is supported on $Q_{k-1,\varepsilon} = \{z | \inf_{v\in Q_{k-1}} |z - v| \le \varepsilon\}$. Since $Q_{k-1,\varepsilon} \cap Q_k \subset B((Mk,0), 16\varepsilon)$, if we choose $k_0$ so that $32\varepsilon \ge 2^{-k_0} \ge 16\varepsilon$, we have that the right-hand side of (5.32) is bounded by

$$
\begin{aligned}
&\sum_{k=k_0}^{\infty} \int_{B(x,2^{-k})\setminus B(x,2^{-k-1})} u(y-x)\, d\mu_{1,\rho}^{X'}(y) \\
&\quad \le c_{24} \sum_{k=k_0}^{\infty} (2^{-k})^{\beta-2} \mu_{1,\rho}^{X'}(B(x,2^{-k})) \\
&\quad \le c_{25} \sum_{k=k_0}^{\infty} 2^{-k(\beta-2)} (2^{-k})^{\beta-\delta} \\
&\quad = c_{25} \sum_{k=k_0}^{\infty} 2^{-k(2\beta-2-\delta)} \le c_{26} \varepsilon^{2\beta-2-\delta}.
\end{aligned}
\tag{5.34}
$$

This completes the proof of Lemma 6. □

If $\varepsilon$ is small enough, we thus conclude from (5.28) and (5.29) that

$$
P(L_k|\mathcal{F}_k) \ge c_{27}\varepsilon^{10+2\beta}\kappa_2/8 \qquad \text{on } \bigcap_{j=1}^{k-1} L_j. \tag{5.35}
$$

Take $\varepsilon$ sufficiently small, but now fix it, and let $\kappa_3 = c_{27}\varepsilon^{4+\beta}\kappa_2/8$. We have

$$
P\Bigg(\bigcap_{j=1}^{k} L_j\Bigg) = E\Bigg[P(L_k|\mathcal{F}_k); \bigcap_{j=1}^{k-1} L_j\Bigg] \ge \kappa_3 P\Bigg(\bigcap_{j=1}^{k-1} L_j\Bigg).
$$

By induction,

$$
P\Bigg(\bigcap_{j=1}^{n} L_j\Bigg) \ge \kappa_3^n.
$$

On the event $M_n = \bigcap_{j=1}^n L_j$, we have that $X_s \in Q_k$ if $k \le s \le k+1$, and so there are no intersections between $X(I_i)$ and $X(I_j)$ if $|i-j| > 1$, where $I_i = [i, i+1]$. Furthermore, on $M_n$, we have

$$
\sum_{k=0}^{n} -\Gamma(I_k) \ge \kappa_1 n/4,
$$

while

$$
\sum_{k=0}^{n} A(I_k; I_{k+1}) \le \kappa_1 n/8.
$$

Since

$$
-\Gamma([0,n]) \ge \sum_{k=0}^{n} -\Gamma(I_k) - \sum_{k=0}^{n} A(I_k; I_{k+1}) \ge \kappa_1 n/8
$$

on the event $M_n$ and $P(M_n) \ge \kappa_3^n$, Theorem 8 is proved. □

**6. The lower tail; $\beta = d$.** In this section we prove Theorem 2 in the critical cases where $\beta = d$. This includes planar Brownian motion and the one-dimensional symmetric Cauchy process.

By the last two lines of Theorem 6, we have

$$E(\alpha(s,t)) = p_1(0)\{(s+t)\log(s+t) - s\log s - t\log t\}. \tag{6.1}$$

Write

$$\eta_t = -\gamma_t - p_1(0)t\log t. \tag{6.2}$$

We have that $\eta_0 = 0$ and, as in the proof of (5.3), for any $s,t>0$, $\eta_{s+t} \le \eta_s + \eta_{s,t}$, where $\eta_{s,t} = -\gamma(\{(u,v)|s \le u \le v \le s+t\}) - p_1(0)t\log t$. For each fixed $s > 0$, $\{\eta_{s,v}; v \ge 0\}$ is independent of $\{\eta_u; u \le s\}$ and $\eta_{s,t} \overset{d}{=} \eta_t$. So by the argument used to obtain (5.9) and (5.10), we obtain

$$E\left(\exp\left\{c\sup_{t\le 1}\eta_t\right\}\right) < \infty \qquad \forall c > 0, \tag{6.3}$$

and

$$\begin{aligned} &E\left(\exp\left\{\frac{1}{p_1(0)}\eta_{s+t}\right\}\right) \\ &\quad\le E\left(\exp\left\{\frac{1}{p_1(0)}\eta_s\right\}\right)E\left(\exp\left\{\frac{1}{p_1(0)}\eta_t\right\}\right) \qquad \forall s,t \ge 0. \end{aligned} \tag{6.4}$$

Therefore, there is a constant $-\infty \le A < \infty$ such that

$$\lim_{t\to\infty} t^{-1}\log E\left(\exp\left\{\frac{1}{p_1(0)}\eta_t\right\}\right) = A \tag{6.5}$$

or, equivalently,

$$\lim_{t\to\infty} t^{-1}\log\left(t^{-t}E\left(\exp\left\{-\frac{1}{p_1(0)}\gamma_t\right\}\right)\right) = A. \tag{6.6}$$

Take $t = n$ to be an integer. By scaling and Stirling's formula,

$$\lim_{n\to\infty}\frac{1}{n}\log\left((n!)^{-1}E\left(\exp\left\{-\frac{n}{p_1(0)}\gamma_1\right\}\right)\right) = A + 1. \tag{6.7}$$

By [12], Lemma 2.3,

$$\lim_{t\to\infty} t^{-1}\log P\left(\exp\left\{-\frac{1}{p_1(0)}\gamma_1\right\} \ge t\right) = -e^{-A-1} \equiv -b_\psi \tag{6.8}$$

or, equivalently,

$$\lim_{t\to\infty} t^{-1}\log P(-\gamma_1 \ge p_1(0)\log t) = -L, \tag{6.9}$$

which proves (1.9). It remains to show that $b_\psi < \infty$. That $b_\psi < \infty$ for the $\beta = d = 2$ case was shown in [2], Section 5. A very similar proof takes care of the $\beta = d = 1$ case. Note that the proof in [2] does not rely on the continuity of Brownian paths. Instead of the $t^{1/2}$ scaling there, we now have $t^1$ scaling. Instead of $1/(2\pi)$, we now have $p_1(0)$, which in the $\beta = d = 1$ case is equal to $1/\pi$. This completes the proof of Theorem 2. □

## 7. The lim sup result.

PROOF OF THEOREM 3. We begin with a lemma.

LEMMA 7. *If $a < a_\psi$, there exists $C < \infty$ such that*

$$P\Big(\sup_{t\le 1} \gamma_t \ge u^{d/\beta}\Big) \le Ce^{-au}, \qquad u > 0. \tag{7.1}$$

PROOF. It follows from (4.8) and scaling that

$$\sup_{t\le 1} P(\gamma_t \ge u^{d/\beta}) \le Ce^{-au}, \qquad u > 0. \tag{7.2}$$

Let $\Gamma([s,t]) := \gamma(\{(u,v) | s \le u \le v \le t\})$. For any $s < t$,

$$\gamma_t - \gamma_s = \gamma([0,s;s,t]) + \Gamma([s,t]), \tag{7.3}$$

with $\gamma([0,s;s,t]) \overset{d}{=} \{\alpha_{s,t-s}\}_0$ and $\Gamma([s,t]) \overset{d}{=} \gamma_{t-s}$.

Using (7.3), it then follows from (2.16) and (3.21) that for some $\theta > 0$,

$$\sup_{s<t\le 1} E\bigg(\exp\bigg\{\theta\bigg|\frac{\gamma_t - \gamma_s}{(t-s)^{1-d/2\beta}}\bigg|^{\beta/d}\bigg\}\bigg) < \infty, \tag{7.4}$$

hence, by Chebyshev, that for some $c > 0$,

$$P(|\gamma_t - \gamma_s| \ge u^{d/\beta}) \le Ce^{-cu/(t-s)^\zeta}, \qquad u > 0, \tag{7.5}$$

uniformly in $0 \le s < t \le 1$, where $\zeta = \beta/d - 1/2 > 0$. Lemma 7 then follows from the chaining argument used in the proof of Proposition 4.1 of [2]. □

It is now straightforward to use scaling and Borel–Cantelli to get the following:

LEMMA 8.

$$\limsup_{t\to\infty} \frac{\gamma_t}{t^{(2-d/\beta)}(\log\log t)^{d/\beta}} \le a_\psi^{-d/\beta} \qquad a.s. \tag{7.6}$$

PROOF. Let $M > 1/a_\psi$. Choose $\varepsilon > 0$ small and $q > 1$ close to 1 so that $M(a_\psi - 2\varepsilon)/q^{2\zeta} > 1$. Let $t_n = q^n$ and let

$$C_n = \Big\{\sup_{s\le t_n} \gamma_s > t_{n-1}^{(2-d/\beta)}(M\log\log t_{n-1})^{d/\beta}\Big\}. \tag{7.7}$$

By Lemma 7 and scaling, the probability of $C_n$ is bounded by

$$c_1 e^{-(a_\psi - \varepsilon)M(t_{n-1}/t_n)^{2\zeta} \log\log t_{n-1}}.$$

By our choices of $\varepsilon$ and $q$, this is summable, so by Borel–Cantelli the probability that $C_n$ happens infinitely often is zero. To complete the proof, we point out that if $\gamma_t > t^{(2-d/\beta)}(M \log\log t)^{d/\beta}$ for some $t \in [t_{n-1}, t_n]$, then the event $C_n$ occurs. This completes the proof of Lemma 8. $\square$

To finish the proof of Theorem 3 we prove the following:

Lemma 9.

$$\limsup_{t\to\infty} \frac{\gamma_t}{t^{(2-d/\beta)}(\log\log t)^{d/\beta}} \ge a_\psi^{-d/\beta} \qquad a.s. \tag{7.8}$$

Proof. Let $a > a_\psi$ and let $a'$ be the midpoint of $(a_\psi, a)$. Then by (4.8),

$$P(\gamma_1 \ge (u \log\log t)^{d/\beta}) \ge c_2 e^{-a' u \log\log t}, \qquad u > 0. \tag{7.9}$$

Let $\delta > 0$ be small enough so that $(1+\delta)a'/a < 1$ and set $t_n = e^{n^{1+\delta}}$. Recall that $\Gamma([s,t]) \stackrel{d}{=} \gamma_{t-s}$. Using (7.9) and scaling, it is straightforward to obtain

$$\sum_{n=1}^{\infty} P\left(\Gamma([t_{n-1}, t_n]) > t_n^{(2-d/\beta)} \left(\frac{\log\log t_n}{a}\right)^{d/\beta}\right) = \infty.$$

Using the fact that different pieces of the path of a stable process are independent and Borel–Cantelli,

$$\limsup_{n\to\infty} \frac{\Gamma([t_{n-1}, t_n])}{t_n^{(2-d/\beta)}(\log\log t_n)^{d/\beta}} > \frac{1}{a^{d/\beta}} \qquad \text{a.s.} \tag{7.10}$$

Let $\varepsilon > 0$. From (3.21), scaling and Borel–Cantelli, it follows that

$$|\Gamma([0, t_{n-1}])| = |\gamma_{t_{n-1}}| = O(\varepsilon t_n^{(2-d/\beta)}(\log\log t_n)^{d/\beta}) \qquad \text{a.s.} \tag{7.11}$$

Since

$$\begin{aligned} \gamma_{t_n} &= \Gamma([0, t_n]) \\ &= \Gamma([t_{n-1}, t_n]) + \Gamma([0, t_{n-1}]) + \gamma([0, t_{n-1}]; [t_{n-1}, t_n]) \end{aligned} \tag{7.12}$$

and $\gamma([0,s];[s,t]) \stackrel{d}{=} \{\alpha_{s,t-s}\}_0$ with $\alpha_{s,t-s} \ge 0$, we have our result from (7.10), (7.11), (7.12) and the fact, from Theorem 6, that

$$E\alpha_{t_{n-1}, t_n - t_{n-1}} \le E\alpha_{t_n} = c_6 t_n^{(2-d/\beta)} = o(t_n^{(2-d/\beta)}(\log\log t_n)^{d/\beta}).$$

This completes the proof of Lemma 9. $\square$

Lemmas 8 and 9 together imply Theorem 3. $\square$

## 8. The lim inf result.

PROOF OF THEOREM 4. We consider first the case when $\beta < d$. Let $D_t = -\gamma_t$. We begin with a lemma.

LEMMA 10. *If $b < b_\psi$, there exists $C < \infty$ such that*

$$P\Big(\sup_{t\le 1} D_t \ge u^{d/\beta-1}\Big) \le Ce^{-bu}, \qquad u > 0. \tag{8.1}$$

PROOF. It follows from (1.8) and scaling (1.10) that

$$\lim_{u\to\infty} u^{-1}\log P(D_1 \ge u^{d/\beta-1}) = -b_\psi. \tag{8.2}$$

Scaling once more shows that, for any $t > 0$,

$$P(D_t \ge u^{d/\beta-1}) \le Ce^{-bu/t^\eta}, \qquad u > 0, \tag{8.3}$$

with $\eta = (2-d/\beta)/(d/\beta-1) > 0$. For any $s < t$,

$$\begin{aligned} D_t - D_s &= -\gamma([0,s;s,t]) - \Gamma([s,t]) \\ &\le E(\alpha_{s,t-s}) - \Gamma([s,t]) \\ &\le c_\psi(t-s)^{2-2/\beta} - \Gamma([s,t]), \end{aligned} \tag{8.4}$$

with $-\Gamma([s,t]) := D_{t-s}$ and we have used Theorem 6

$$E(\alpha_{s,t-s}) = c_\psi[s^{2-2/\beta} + (t-s)^{2-2/\beta} - t^{2-2/\beta}] \le c_\psi(t-s)^{2-2/\beta}. \tag{8.5}$$

Lemma 10 then follows from the chaining argument used in the proof of Proposition 4.1 of [2]. □

It is now straightforward to use scaling and Borel–Cantelli to get the following:

LEMMA 11.

$$\limsup_{t\to\infty} \frac{D_t}{t^{(2-d/\beta)}(\log\log t)^{d/\beta-1}} \le b_\psi^{-(d/\beta-1)} \qquad a.s. \tag{8.6}$$

PROOF. Let $M > 1/b_\psi$. Choose $\varepsilon > 0$ small and $q > 1$ close to 1 so that $M(b_\psi - 2\varepsilon)/q^\rho > 1$. Let $t_n = q^n$ and let

$$C_n = \Big\{\sup_{s\le t_n} D_s > t_{n-1}^{(2-d/\beta)}(M\log\log t_{n-1})^{d/\beta-1}\Big\}. \tag{8.7}$$

By Lemma 7 and scaling, the probability of $C_n$ is bounded by

$$c_1 e^{-(b_\psi-\varepsilon)M(t_{n-1}/t_n)^\rho \log\log t_{n-1}}.$$

By our choices of $\varepsilon$ and $q$, this is summable, so by Borel–Cantelli the probability that $C_n$ happens infinitely often is zero. To complete the proof, we point out that if $D_t > t^{(2-d/\beta)}(M\log\log t)^{d/\beta-1}$ for some $t \in [t_{n-1}, t_n]$, then the event $C_n$ occurs. This completes the proof of Lemma 11. $\square$

To finish the proof of Theorem 4 when $\beta < d$, we prove the next lemma.

LEMMA 12.

$$\limsup_{t\to\infty} \frac{D_t}{t^{(2-d/\beta)}(\log\log t)^{d/\beta-1}} \geq b_\psi^{-(d/\beta-1)} \qquad a.s. \tag{8.8}$$

PROOF. Let $b > b_\psi$ and let $b'$ be the midpoint of $(b_\psi, b)$. Then by (8.2),

$$P(D_1 \geq (u\log\log t)^{d/\beta-1}) \geq c_2 e^{-b'u\log\log t}, \qquad u > 0. \tag{8.9}$$

Let $\delta > 0$ be small enough so that $(1+\delta)b'/b < 1$ and set $t_n = e^{n^{1+\delta}}$. Recall that $\Gamma([s,t]) \stackrel{d}{=} \gamma_{t-s}$. Using (8.9) and scaling, it is straightforward to obtain

$$\sum_{n=1}^{\infty} P\biggl(-\Gamma([t_{n-1}, t_n]) > t_n^{(2-d/\beta)}\biggl(\frac{\log\log t_n}{b}\biggr)^{d/\beta-1}\biggr) = \infty.$$

Using the fact that different pieces of the path of a stable process are independent and Borel–Cantelli,

$$\limsup_{n\to\infty} \frac{-\Gamma([t_{n-1}, t_n])}{t_n^{(2-d/\beta)}(\log\log t_n)^{d/\beta-1}} > \frac{1}{b^{d/\beta-1}} \qquad \text{a.s.} \tag{8.10}$$

Let $\varepsilon > 0$. From (3.21), scaling and Borel–Cantelli, it follows that

$$|\Gamma([0, t_{n-1}])| = |\gamma_{t_{n-1}}| = O(\varepsilon t_n^{(2-d/\beta)}(\log\log t_n)^{d/\beta-1}) \qquad \text{a.s.} \tag{8.11}$$

Note that

$$\begin{aligned} D_{t_n} &= -\Gamma([0, t_n]) \\ &= -\Gamma([t_{n-1}, t_n]) - \Gamma([0, t_{n-1}]) - \gamma([0, t_{n-1}]; [t_{n-1}, t_n]) \end{aligned} \tag{8.12}$$

and $\gamma([0,s];[s,t]) \stackrel{d}{=} \{\alpha_{s,t-s}\}_0$. Using (2.16),

$$\begin{aligned} &P(\alpha([0, t_{n-1}]; [t_{n-1}, t_n]) > t_n^{(2-d/\beta)}) \\ &\quad \leq P\biggl(\frac{\alpha([0, t_{n-1}]; [t_{n-1}, t_n])}{(t_{n-1}(t_n - t_{n-1}))^{(1-d/2\beta)}} > (t_n/t_{n-1})^{(1-d/2\beta)}\biggr) \\ &\quad \leq e^{-(t_n/t_{n-1})^{(\beta/d-1/2)}}, \end{aligned} \tag{8.13}$$

which is summable. Using Borel–Cantelli, we have

$$\alpha([0, t_{n-1}]; [t_{n-1}, t_n]) = o(t_n^{(2-d/\beta)}(\log\log t_n)^{d/\beta-1}). \tag{8.14}$$

Substituting this, (8.10) and (8.11) in (8.12) completes the proof of Lemma 12.

□

Lemmas 11 and 12 together imply Theorem 4 when $\beta < d$. The case of $\beta = d$ follows from (6.9) and the proof of [2], Theorem 1.5. □

**Acknowledgment.** We thank Evarist Giné for supplying the elegant proof of Lemma 1.

R. Bass
Department of Mathematics
University of Connecticut
Storrs, Connecticut 06269-3009
USA
e-mail: bass@math.uconn.edu

X. Chen
Department of Mathematics
University of Tennessee
Knoxville, Tennessee 37996-1300
USA
e-mail: xchen@math.utk.edu

J. Rosen
Department of Mathematics
College of Staten Island, CUNY
Staten Island, New York 10314
USA
e-mail: jrosen3@earthlink.net